\documentclass[a4paper,11pt]{article}
\usepackage{graphicx}
\usepackage{amsmath}
\usepackage{amsfonts}
\usepackage{amssymb}
\newtheorem{theorem}{Theorem}

\newtheorem{definition}[theorem]{Definition}

\newtheorem{lemma}[theorem]{Lemma}
\newtheorem{notation}[theorem]{Notation}

\newtheorem{proposition}[theorem]{Proposition}
\newtheorem{remark}[theorem]{Remark}

\newenvironment{proof}[1][Proof]{\textbf{#1.} }{\ \rule{0.5em}{0.5em}}
\begin{document}

\title{Remarks on the embedding of spaces of distributions into spaces of Colombeau
generalized functions}
\author{A. Delcroix\\Equipe Analyse Alg\'{e}brique Non Lin\'{e}aire\\\textit{Laboratoire Analyse, Optimisation, Contr\^{o}le}\\Facult\'{e} des sciences - Universit\'{e} des Antilles et de la Guyane\\97159\ Pointe-\`{a}-Pitre Cedex Guadeloupe}
\maketitle
\begin{abstract}
We present some remarks about the embedding of spaces of Schwartz
distributions into spaces of Colombeau generalized functions. Following ideas
of M.\ Nedeljkov \emph{et alii}, we recall how a good choice of compactly
supported mollifiers allows to perform globally the embedding \ of
$\mathcal{D}^{\prime}\left(  \Omega\right)  $ into $\mathcal{G}\left(
\Omega\right)  $.\ We show that this embedding is equal to the one obtained
with local and sheaf arguments by M.~Grosser \textit{et alii}, this giving
various equivalent technics to embed $\mathcal{D}^{\prime}\left(
\Omega\right)  $ into $\mathcal{G}\left(  \Omega\right)  .$
\end{abstract}

\noindent\textbf{Mathematics Subject Classification (2000): 46E10, 46E25,
46F05, 46F30\smallskip}

\bigskip

\section{Introduction\label{SPInt}}

The question of embedding classical spaces such as $\mathrm{C}^{0}\left(
\Omega\right)  $, $\mathrm{C}^{\infty}\left(  \Omega\right)  $, $\mathcal{D}%
^{\prime}\left(  \Omega\right)  $ (where $\Omega$ is an open subset of
$\mathbb{R}^{d}$, $d\in\mathbb{N}$) in spaces of generalized functions arises
naturally. \smallskip

The embedding of $\mathrm{C}^{\infty}\left(  \Omega\right)  $ into
$\mathcal{G}\left(  \Omega\right)  $ is classically done by the canonical map
\[
\sigma:\mathrm{C}^{\infty}\left(  \Omega\right)  \rightarrow\mathcal{G}\left(
\Omega\right)  \;\;\ f\rightarrow\left[  \left(  f_{\varepsilon}\right)
_{\varepsilon}\right]  \;\mathrm{with\;}f_{\varepsilon}=f\text{ for
}\varepsilon\in\left(  0,1\right]
\]
which is an injective homomorphism of algebras. ($\left[  \left(
f_{\varepsilon}\right)  _{\varepsilon}\right]  $ denotes the class of $\left(
f_{\varepsilon}\right)  _{\varepsilon}$ in the factor algebra $\mathcal{G}%
\left(  \Omega\right)  $: See section \ref{SPPrel} for a short presentation of
$\mathcal{G}\left(  \Omega\right)  $ or \cite{Col1}, \cite{GKOS} for a
complete construction.)\smallskip

For the embedding of $\mathcal{D}^{\prime}\left(  \Omega\right)  $ into
$\mathcal{G}\left(  \Omega\right)  $ the following additional assumption is
required: If $\iota_{A}$ is the expected embedding, one wants the following
diagram to be commutative
\begin{equation}%
\begin{array}
[c]{ccc}%
\mathrm{C}^{\infty}\left(  \Omega\right)  & \longrightarrow & \mathcal{D}%
^{\prime}\left(  \Omega\right) \\
& \searrow\!\sigma & \downarrow\iota_{A}\\
&  & \mathcal{G}\left(  \Omega\right)
\end{array}
\label{Scomdiag}%
\end{equation}
that is $\iota_{A\left|  \mathrm{C}^{\infty}\left(  \Omega\right)  \right.
}=\sigma$.\smallskip

In \cite{GKOS}, this program is fulfilled by using the sheaf properties of
Colombeau simplified algebras. Let us quote the main step of the
construction.\ First, an embedding $\iota_{0}$ of $\mathcal{E}^{\prime}\left(
\mathbb{R}^{d}\right)  $ into $\mathcal{G}\left(  \mathbb{R}^{d}\right)  $ is
realized by convolution of compactly supported distributions with suitable
mollifiers $\left(  \rho_{\varepsilon}\right)  _{\varepsilon}$ belonging to
$\mathcal{S}\left(  \mathbb{R}^{d}\right)  $. In fact, this map $\iota_{0}$
can be considered as an embedding of $\mathcal{E}^{\prime}\left(
\mathbb{R}^{d}\right)  $ into $\mathcal{G}_{C}\left(  \mathbb{R}^{d}\right)
$, the subalgebra of $\mathcal{G}\left(  \mathbb{R}^{d}\right)  $ of compactly
supported generalized functions, since the support of $T\in\mathcal{E}%
^{\prime}\left(  \mathbb{R}^{d}\right)  $ is equal to the support of its image
by $\iota_{0}$ (Proposition 1.2.12\ of \cite{GKOS}). The following step of the
construction of $\iota_{A}$ is to consider for every open set $\Omega
\subset\mathbb{R}^{d}$ an open covering $\left(  \Omega_{\lambda}\right)
_{\lambda}$ of $\Omega$ with relatively compact open sets and to embed
$\mathcal{D}^{\prime}\left(  \Omega\right)  $ into $\mathcal{G}\left(
\Omega_{\lambda}\right)  $ with the help of cutoff functions and $\iota_{0}$.
Using a partition of unity subordinate to $\left(  \Omega_{\lambda}\right)
_{\lambda}$, $\iota_{A}$ is constructed by ``gluing the bits obtained before
together''. Finally, it is shown that the embedding $\iota_{A}$ does not
depend on the choice of $\left(  \Omega_{\lambda}\right)  _{\lambda}$ and
other material of the construction, excepted the net $\left(  \rho
_{\varepsilon}\right)  _{\varepsilon}$.\smallskip

On one hand, the choice of \textit{ad hoc} not compactly supported mollifiers
renders trivially the diagram (\ref{Scomdiag}) commutative. On the other hand,
$\rho_{\varepsilon}$ cannot be convoluted with elements of $\mathcal{D}%
^{\prime}\left(  \Omega\right)  $ unrestrictedly, obliging to consider first
compactly supported distributions, and then sheaf arguments.\smallskip

In \cite{NePiSc}, the authors give an other construction which avoid this
drawback. The main idea is to use compactly supported mollifiers enough close
from the \textit{ad hoc} mollifiers $\left(  \rho_{\varepsilon}\right)
_{\varepsilon}$ of \cite{GKOS}.\ This is done by regular cutoff of\ $\left(
\rho_{\varepsilon}\right)  _{\varepsilon}$, this cutoff being defined with an
other rate of growth than the net $\left(  \rho_{\varepsilon}\right)
_{\varepsilon}$, let say in $\ln\varepsilon$ whereas the scale of growth of
$\left(  \rho_{\varepsilon}\right)  _{\varepsilon}$ is in $1/\varepsilon$.
This permits to keep the good properties of the embedding in particular the
commutativity of the diagram (\ref{Scomdiag}). We present in details this
construction in section \ref{SPEmbR} for the case $\Omega=\mathbb{R}^{d}$.\smallskip

In section \ref{SPComp}, we show that these embeddings are in fact equal,
consequently only depending on the choice of the mollifiers $\left(
\rho_{\varepsilon}\right)  _{\varepsilon}$. (This dependance is well known for
the simplified Colombeau algebra.)\ We finally turn to the case of the
embedding of $\mathcal{D}^{\prime}\left(  \Omega\right)  $ into the simplified
Colombeau Algebra $\mathcal{G}\left(  \Omega\right)  $ where $\Omega$ is an
arbitrary open subset of $\mathbb{R}^{d}$ (Section \ref{SPEmbOm}). We show
that for the global construction of \cite{NePiSc} an additional cutoff applied
to the elements of $\mathcal{D}^{\prime}\left(  \Omega\right)  $ is needed. We
also give a local version (with no cutoff on the distribution) of the
construction of \cite{NePiSc}.\smallskip

This note comes from a workshop in Paris 7 and seminars of the team AANL of
the laboratory AOC held in June and September 2003. I deeply think
D.\ Scarpalezos\ and J.-A.\ Marti for the discussions about these constructions.

\section{Preliminaries\label{SPPrel}}

\subsection{The sheaf of Colombeau simplified algebras\label{SBColo}}

Let $\mathrm{C}^{\infty}$ be the sheaf of complex valued smooth functions on
$\mathbb{R}^{d}$ ($d\in\mathbb{N}$) with the usual topology of uniform
convergence.\ For every open set $\Omega$ of $\mathbb{R}^{d}$, this topology
can be described by the family of semi norms
\[
p_{K,l}(f)=\sup_{\left|  \alpha\right|  \leq l,K\Subset\Omega}\left|
\partial^{\alpha}f\left(  x\right)  \right|  ,
\]
where the notation $K\Subset\Omega$ means that\ the set $K$ is a compact set
included in $\Omega$.

Let us set
\begin{align*}
&  \!\!\!\!\!\!\!\!\!\!\mathcal{F}\left(  \mathrm{C}^{\infty}\left(
\Omega\right)  \right) \\
~~~~  &  =\left\{  \left(  f_{\varepsilon}\right)  _{\varepsilon}\in
\mathrm{C}^{\infty}\left(  \Omega\right)  ^{\left(  0,1\right]  }\,\left|
\,\forall l\in\mathbb{N},\;\forall K\Subset\Omega,\;\exists q\in
\mathbb{N},\;p_{K,l}\left(  f_{\varepsilon}\right)  =\mathrm{o}\left(
\varepsilon^{-q}\right)  \;\mathrm{for}\;\varepsilon\rightarrow0\right.
\right\}  ,\\
&  \!\!\!\!\!\!\!\!\!\!\mathcal{N}\left(  \mathrm{C}^{\infty}\left(
\Omega\right)  \right) \\
~~~~  &  =\left\{  \left(  f_{\varepsilon}\right)  _{\varepsilon}\in
\mathrm{C}^{\infty}\left(  \Omega\right)  ^{\left(  0,1\right]  }\,\left|
\,\forall l\in\mathbb{N},\;\forall K\Subset\Omega,\;\forall p\in
\mathbb{N},\;p_{K,l}\left(  f_{\varepsilon}\right)  =\mathrm{o}\left(
\varepsilon^{p}\right)  \;\mathrm{for}\;\varepsilon\rightarrow0\right.
\right\}  .
\end{align*}

\begin{lemma}
\cite{JAM1} and \cite{JAM2}\newline $i$.~The functor $\mathcal{F}%
:\Omega\rightarrow\mathcal{F}\left(  \mathrm{C}^{\infty}\left(  \Omega\right)
\right)  $ defines a sheaf of subalgebras of the sheaf $\left(  \mathrm{C}%
^{\infty}\right)  ^{\left(  0,1\right]  }$\newline $ii$.~The functor
$\mathcal{N}:\Omega\rightarrow\mathcal{N}\left(  \mathrm{C}^{\infty}\left(
\Omega\right)  \right)  $ defines a sheaf of ideals of the sheaf $\mathcal{F}$.
\end{lemma}

We shall note prove in detail this lemma but quote the two mains arguments:

\noindent$i$.~for each open subset $\Omega$ of $X$, the family of seminorms
$\left(  p_{K,l}\right)  $ related to $\Omega$ is compatible with the
algebraic structure of $\mathrm{C}^{\infty}\left(  \Omega\right)  ;$\ In
particular:
\[
\forall l\in\mathbb{N},\;\forall K\Subset\Omega,\;\;\exists C\in\mathbb{R}%
_{+}^{\ast},\;\;\forall\left(  f,g\right)  \in\left(  \mathrm{C}^{\infty
}\left(  \Omega\right)  \right)  ^{2}\;\;p_{K,l}\left(  fg\right)  \leq
Cp_{K,l}\left(  f\right)  p_{K,l}\left(  g\right)  ,
\]

\noindent$ii$.~For two open subsets $\Omega_{1}\subset\Omega_{2}$ of
$\mathbb{R}^{d}$, the family of seminorms $\left(  p_{K,l}\right)  $ related
to $\Omega_{1}$ is included in the family of seminorms related to $\Omega_{2}$
and
\[
\forall l\in\mathbb{N},\;\forall K\Subset\Omega_{1},\;\;\forall f\in
\mathrm{C}^{\infty}\left(  \Omega_{2}\right)  ,\;\;p_{K,l}\left(  f_{\left|
\Omega_{1}\right.  }\right)  =p_{K,l}\left(  f\right)  .
\]

\begin{definition}
The sheaf of factor algebras
\[
\mathcal{G}\left(  \mathrm{C}^{\infty}\left(  \cdot\right)  \right)
=\mathcal{F}\left(  \mathrm{C}^{\infty}\left(  \cdot\right)  \right)
/\mathcal{N}\left(  \mathrm{C}^{\infty}\left(  \cdot\right)  \right)
\]
is called the sheaf of \emph{Colombeau simplified algebras}.
\end{definition}

The sheaf $\mathcal{G}$ turns to be a sheaf of differential algebras and a
sheaf of modulus on the factor ring $\overline{\mathbb{C}}=\mathcal{F}\left(
\mathbb{C}\right)  /\mathcal{N}\left(  \mathbb{C}\right)  $ with
\begin{align*}
\mathcal{F}\left(  \mathbb{K}\right)   &  =\left\{  \left(  r_{\varepsilon
}\right)  _{\varepsilon}\in\mathbb{K}^{\left(  0,1\right]  }\,\left|
\,\exists q\in\mathbb{N},\;\left|  r_{\varepsilon}\right|  =\mathrm{o}\left(
\varepsilon^{-q}\right)  \;\mathrm{for}\;\varepsilon\rightarrow0\right.
\right\}  ,\\
\mathcal{N}\left(  \mathbb{K}\right)   &  =\left\{  \left(  r_{\varepsilon
}\right)  _{\varepsilon}\in\mathbb{K}^{\left(  0,1\right]  }\,\left|
\,\forall p\in\mathbb{N},\;\left|  r_{\varepsilon}\right|  =\mathrm{o}\left(
\varepsilon^{p}\right)  \;\mathrm{for}\;\varepsilon\rightarrow0\right.
\right\}  ,
\end{align*}
with $\mathbb{K}=\mathbb{C}$ or $\mathbb{K}=\mathbb{R},\;\mathbb{R}_{+}$.

\begin{notation}
In the sequel we shall note, as usual, $\mathcal{G}\left(  \Omega\right)  $
instead of $\mathcal{G}\left(  \mathrm{C}^{\infty}\left(  \Omega\right)
\right)  $ the algebra of generalized functions on $\Omega$.
\end{notation}

\subsection{Local structure of distributions\label{SBLocDP}}

To fix notations, we recall here two classical results on the local structure
of distributions which are going to be used in the sequel. We refer the reader
to \cite{Schwartz1} chapter 3, specially theorem XXI \& XXVI, \ for proofs and
details. Let $\Omega$ be an open subset of $\mathbb{R}^{d}$ ($d\in\mathbb{N}$).

\begin{theorem}
\label{SlocStrucDP}For all $T\in\mathcal{D}^{\prime}\left(  \Omega\right)  $
and all $\Omega^{\prime}$ open subset of $\mathbb{R}^{d}$ with $\overline
{\Omega^{\prime}}\Subset\Omega$, there exists $f\in\mathrm{C}^{0}\left(
\mathbb{R}^{d}\right)  $ whose support is contained in an arbitrary
neighborhood of $\overline{\Omega^{\prime}}$, $\alpha\in\mathbb{N}^{d}$ such
that $T_{\left|  \Omega^{\prime}\right.  }=\partial^{\alpha}f.$
\end{theorem}

\begin{theorem}
\label{SlocStrucEP}For all $T\in\mathcal{E}^{\prime}\left(  \Omega\right)  $,
there exists an integer $r\geq0$, a finite family $\left(  f_{\alpha}\right)
_{0\leq\left|  \alpha\right|  \leq r}$ ($\alpha\in\mathbb{N}^{d}$) with each
$f_{\alpha}\in\mathrm{C}^{0}\left(  \mathbb{R}^{d}\right)  $ having its
support contained in the same arbitrary neighborhood of the support of $T$,
such that $T=\sum_{0\leq\left|  \alpha\right|  \leq r}\partial^{\alpha
}f_{\alpha}.$
\end{theorem}

\section{Embedding of $\mathcal{D}^{\prime}\left(  \mathbb{R}^{d}\right)  $ in
$\mathcal{G}\left(  \mathbb{R}^{d}\right)  $\label{SPEmbR}}

\subsection{Construction of the mollifiers\label{SBMoll}}

Take $\rho\in\mathcal{S}\left(  \mathbb{R}^{d}\right)  $ even such that
\begin{equation}
\int\rho\left(  x\right)  \,\mathrm{d}x=1,\ \ \ \ \ \ \ \ \ \ \ \int x^{m}%
\rho\left(  x\right)  \,\mathrm{d}x=0\;\mathrm{for\;all\;}m\in\mathbb{N}%
^{d}\backslash\left\{  0\right\}  \label{SBinfini}%
\end{equation}
and $\chi\in\mathcal{D}\left(  \mathbb{R}^{d}\right)  $ such that $0\leq
\chi\leq1$,$\;\chi\equiv1$ on $\overline{B(0,1)}$ and $\chi\equiv0$ on
$\mathbb{R}^{d}\backslash B(0,2)$.\ Define
\[
\forall\varepsilon\in\left(  0,1\right]  ,\;\;\forall x\in\mathbb{R}%
^{d},\;\;\;\;\rho_{\varepsilon}\left(  x\right)  =\frac{1}{\varepsilon^{d}%
}\rho\left(  \frac{x}{\varepsilon}\right)
\]
and
\[
\forall\varepsilon\in\left(  0,1\right)  ,\;\;\forall x\in\mathbb{R}%
^{d},\;\;\;\;\theta_{\varepsilon}\left(  x\right)  =\rho_{\varepsilon}\left(
x\right)  \chi\left(  x\left|  \ln\varepsilon\right|  \right)
\,;\;\;\text{\ \ }\theta_{1}\left(  x\right)  =1\text{.}%
\]

\begin{remark}
\label{RemarkMod}The nets $\left(  \rho_{\varepsilon}\right)  _{\varepsilon} $
and $\left(  \theta_{\varepsilon}\right)  _{\varepsilon}$ defined above belong
to $\mathcal{F}\left(  \mathrm{C}^{\infty}\left(  \mathbb{R}^{d}\right)
\right)  $.
\end{remark}

Let us verify this result for $\left(  \theta_{\varepsilon}\right)
_{\varepsilon}$ and $d=1$. Fixing $\alpha\in\mathbb{N}$, we have
\begin{align*}
\left.  \forall x\in\mathbb{R},\;\;\;\;\partial^{\alpha}\theta_{\varepsilon
}\left(  x\right)  \right.   &  =\sum_{\beta=0}^{\alpha}\mathrm{C}_{\alpha
}^{\beta}\partial^{\beta}\rho_{\varepsilon}\left(  x\right)  \partial
^{\alpha-\beta}\chi\left(  x\left|  \ln\varepsilon\right|  \right)  ,\\
&  =\sum_{\beta=0}^{\alpha}\mathrm{C}_{\alpha}^{\beta}\varepsilon^{-1-\beta
}\left|  \ln\varepsilon\right|  ^{\alpha-\beta}\rho^{\left(  \beta\right)
}\left(  \frac{x}{\varepsilon}\right)  \chi^{\left(  \alpha-\beta\right)
}\left(  x\left|  \ln\varepsilon\right|  \right)  .
\end{align*}
For all $\beta\in\left\{  0,\ldots,\alpha\right\}  $, we have $\varepsilon
^{-1-\beta}\left|  \ln\varepsilon\right|  ^{\alpha-\beta}=\mathrm{o}\left(
\varepsilon^{-2-\alpha}\right)  $ for $\varepsilon\rightarrow0$. As
$\rho^{\left(  n\right)  }$ and $\chi^{\left(  n\right)  }$ are bounded, there
exists $C(\alpha)$ such that
\[
\forall x\in\mathbb{R},\;\;\;\;\left|  \partial_{\varepsilon}^{\alpha}\left(
\theta\left(  x\right)  \right)  \right|  \leq C(\alpha)\varepsilon
^{-2-\alpha}.
\]
Our claim follows from this last inequality.

\begin{lemma}
\label{SLmnGoodM}With the previous notations, the following properties holds
\begin{gather}
\left(  \theta_{\varepsilon}\right)  _{\varepsilon}-\left(  \rho_{\varepsilon
}\right)  _{\varepsilon}\in\mathcal{N}\left(  \mathrm{C}^{\infty}\left(
\mathbb{R}^{d}\right)  \right)  ,\label{SLmnGoodM1}\\
\forall k\in\mathbb{N},\;\;\;\int\theta_{\varepsilon}\left(  x\right)
\,\mathrm{d}x=1+\mathrm{o}\left(  \varepsilon^{k}\right)  \;\text{for}%
\;\varepsilon\rightarrow0,\label{SLmnGoodM1B}\\
\forall k\in\mathbb{N},\;\;\forall m\in\mathbb{N}^{d}\backslash\left\{
0\right\}  ,\;\ \ \int x^{m}\theta_{\varepsilon}\left(  x\right)
\,\mathrm{d}x=\mathrm{o}\left(  \varepsilon^{k}\right)  \;\text{for}%
\;\varepsilon\rightarrow0. \label{SLmnGoodM2B}%
\end{gather}
\end{lemma}

In other words, we have
\[
\left(  \int\theta_{\varepsilon}\left(  x\right)  \,\mathrm{d}x-1\right)
_{\varepsilon}\in\mathcal{N}\left(  \mathbb{R}\right)  \;\;\;\;\;\;\forall
m\in\mathbb{N}^{d}\backslash\left\{  0\right\}  ,\;\left(  \int x^{m}%
\theta_{\varepsilon}\left(  x\right)  \,\mathrm{d}x\right)  _{\varepsilon}%
\in\mathcal{N}\left(  \mathbb{R}\right)  .
\]

\begin{proof}
We consider the case $d=1$ in order to simplify notations.

\noindent\emph{First assertion}.- We have, for all $x\in\mathbb{R}$ and
$\varepsilon\in\left(  0,1\right)  $,
\begin{equation}
\left|  \rho_{\varepsilon}\left(  x\right)  -\theta_{\varepsilon}\left(
x\right)  \right|  \leq\frac{1}{\varepsilon}\left|  \rho\left(  \frac
{x}{\varepsilon}\right)  \right|  \left(  1-\chi\left(  x\left|
\ln\varepsilon\right|  \right)  \right)  \,\leq\frac{1}{\varepsilon}\left|
\rho\left(  \frac{x}{\varepsilon}\right)  \right|  , \label{Cestimate0}%
\end{equation}
Since $\rho$ belongs to $\mathcal{S}\left(  \mathbb{R}\right)  $, for all
integer $k>0$ there exists a constant $C(k)$ such that
\[
\forall x\in\mathbb{R},\;\;\;\left|  \rho\left(  x\right)  \right|  \leq
\frac{C(k)}{1+\left|  x\right|  ^{k}}\text{.}%
\]
Then, for all $x\in\mathbb{R}$ with $\left|  x\right|  \geq1/\left|
\ln\varepsilon\right|  ,$%
\begin{equation}
\frac{1}{\varepsilon}\left|  \rho\left(  \frac{x}{\varepsilon}\right)
\right|  \leq\frac{C(k)}{\varepsilon^{k}+\left|  x\right|  ^{k}}%
\varepsilon^{k-1}\leq C(k)\left|  \ln\varepsilon\right|  ^{k}\varepsilon
^{k-1}=\mathrm{o}\left(  \varepsilon^{k-2}\right)  \label{SCestimate}%
\end{equation}
According to remark \ref{RemarkMod}, $\left(  \rho_{\varepsilon}%
-\theta_{\varepsilon}\right)  _{\varepsilon}\in\mathcal{F}\left(
\mathrm{C}^{\infty}\left(  \mathbb{R}\right)  \right)  $.\ Then we can
conclude without estimating the derivatives that $\left(  \rho_{\varepsilon
}-\theta_{\varepsilon}\right)  _{\varepsilon}\in\mathcal{N}\left(
\mathrm{C}^{\infty}\left(  \mathbb{R}\right)  \right)  $ by using theorem
1.2.3. of \cite{GKOS}.\medskip

\noindent\emph{Second and third assertions}.- According definition of $\rho$
we find that
\[
\int\rho_{\varepsilon}\left(  x\right)  \,\mathrm{d}x=1\,;\ \;\;\;\;\forall
m\in\mathbb{N}\backslash\left\{  0\right\}  ,\;\int x^{m}\rho_{\varepsilon
}\left(  x\right)  \,\mathrm{d}x=0.
\]
So, it suffices to show that, for all $m\in\mathbb{N}$ and $k\in\mathbb{N}$,
\[
\Delta_{m,\varepsilon}\left(  x\right)  =\int x^{m}\left(  \rho_{\varepsilon
}\left(  x\right)  -\theta_{\varepsilon}\left(  x\right)  \right)
\,\mathrm{d}x=\mathrm{o}\left(  \varepsilon^{k}\right)  \text{ for
}\varepsilon\rightarrow0
\]
to complete our proof. Let fix $m\in\mathbb{N}$. We have
\[
\Delta_{m,\varepsilon}\left(  x\right)  =\underset{(\ast)}{\underbrace
{\int_{-\infty}^{1/\left|  \ln\varepsilon\right|  }x^{m}\left(  \rho
_{\varepsilon}\left(  x\right)  -\theta_{\varepsilon}\left(  x\right)
\right)  \,\mathrm{d}x}}+\underset{(\ast\ast)}{\underbrace{\int_{1/\left|
\ln\varepsilon\right|  }^{+\infty}x^{m}\left(  \rho_{\varepsilon}\left(
x\right)  -\theta_{\varepsilon}\left(  x\right)  \right)  \,\mathrm{d}x}}.
\]
Let us find an estimate of $(\ast\ast)$. We have, according\ to relations
(\ref{Cestimate0}) and (\ref{SCestimate}),
\[
\forall x>0,\;\;\forall k\in\mathbb{N},\;\;x^{m}\left|  \rho_{\varepsilon
}\left(  x\right)  -\theta_{\varepsilon}\left(  x\right)  \right|  \leq
\frac{1}{\varepsilon}x^{m}\left|  \rho\left(  \frac{x}{\varepsilon}\right)
\right|  \leq C(k)\varepsilon^{k-1}x^{m-k}.
\]
Therefore, by choosing $k\geq m+2$,
\begin{multline*}
\left|  \int_{1/\left|  \ln\varepsilon\right|  }^{+\infty}x^{m}\left(
\rho_{\varepsilon}\left(  x\right)  -\theta_{\varepsilon}\left(  x\right)
\right)  \,\mathrm{d}x\right|  \leq C(k)\varepsilon^{k-1}\int_{1/\left|
\ln\varepsilon\right|  }^{+\infty}x^{m-k}\,\mathrm{d}x\\
\leq\frac{C(k)}{k-1-m}\varepsilon^{k-1}\left|  \ln\varepsilon\right|
^{k-m-1}=\mathrm{o}\left(  \varepsilon^{k-2}\right)  ,\;\text{for }%
\varepsilon\rightarrow0.
\end{multline*}
With a similar estimate for $(\ast)$\ we obtain our claim.
\end{proof}

\subsection{Construction of the embedding $\iota_{S}$\label{SBembR}}

\begin{proposition}
\label{SLmnEmbed}With notations of lemma \ref{SLmnGoodM}, the map
\[
\iota_{S}:\mathcal{D}^{\prime}\left(  \mathbb{R}^{d}\right)  \rightarrow
\mathcal{G}\left(  \mathbb{R}^{d}\right)  \;\;\;T\mapsto\left(  T\ast
\theta_{\varepsilon}\right)  _{\varepsilon}+\mathcal{N}\left(  \mathrm{C}%
^{\infty}\left(  \mathbb{R}^{d}\right)  \right)
\]
is an injective homomorphism of vector spaces.\ Moreover $\iota_{S\left|
\mathrm{C}^{\infty}\left(  \mathbb{R}^{d}\right)  \right.  }=\sigma$.
\end{proposition}

\begin{proof}
We have first to show that for all $T\in\mathcal{D}^{\prime}\left(
\mathbb{R}^{d}\right)  $, $\left(  T\ast\theta_{\varepsilon}\right)
_{\varepsilon}\in\mathcal{F}\left(  \mathrm{C}^{\infty}\left(  \Omega\right)
\right)  $. (This allows to define the map $\iota_{S}$.) Let us fix a compact
$K$.\ Consider $\Omega$ open subset of $\mathbb{R}^{d}$ such that
$K\subset\Omega\subset\overline{\Omega}\Subset\mathbb{R}^{d}$. Let us recall
that
\[
\forall y\in\mathbb{R}^{d},\;\;\;T\ast\theta_{\varepsilon}(y)=\left\langle
T,\left\{  x\mapsto\theta_{\varepsilon}\left(  y-x\right)  \right\}
\right\rangle
\]
For $y\in K$ and $x\in\mathbb{R}^{d}$, we have
\begin{equation}
\theta_{\varepsilon}\left(  y-x\right)  \neq0\Rightarrow y-x\in B(0,\frac
{2}{\left|  \ln\varepsilon\right|  })\Rightarrow x\in B(y,\frac{2}{\left|
\ln\varepsilon\right|  })\Rightarrow x\in\Omega, \label{SsuppTheta}%
\end{equation}
for $\varepsilon$ small enough.

Then, the function $x\mapsto\theta_{\varepsilon}\left(  y-x\right)  $ belongs
to $\mathcal{D}\left(  \Omega\right)  $ and
\[
\left\langle T,\theta_{\varepsilon}\left(  y-\cdot\right)  \right\rangle
=\left\langle T_{\left|  \Omega\right.  },\theta_{\varepsilon}\left(
y-\cdot\right)  \right\rangle .
\]
Using theorem \ref{SlocStrucDP}, we can write $T_{\left|  \Omega\right.
}=\partial_{x}^{\alpha}f$ where $f$ is a continuous function compactly
supported. Then $T\ast\theta_{\varepsilon}=f\ast\partial^{\alpha}%
\theta_{\varepsilon}$ and
\[
\forall y\in K,\;\;\;\left(  T\ast\theta_{\varepsilon}\right)  (y)=\int
_{\Omega}f\left(  y-x\right)  \partial^{\alpha}\theta_{\varepsilon}\left(
x\right)  \,\mathrm{d}x.
\]
According to remark \ref{RemarkMod}, there exists $m\left(  \alpha\right)
\in\mathbb{N}$ such that
\[
\forall x\in\mathbb{R}^{d},\;\;\;\left|  \partial^{\alpha}\theta_{\varepsilon
}\left(  x\right)  \right|  \leq C\varepsilon^{-m\left(  \alpha\right)  }.
\]
We get
\[
\forall y\in K,\;\;\;\left|  \left(  T\ast\theta_{\varepsilon}\right)
(y)\right|  \leq C\sup_{\xi\in\overline{\Omega}}\left|  f\left(  \xi\right)
\right|  \,\mathrm{vol}\left(  \overline{\Omega}\right)  \,\varepsilon
^{-m\left(  \alpha\right)  },
\]
and $\sup_{y\in K}\left|  \left(  T\ast\theta_{\varepsilon}\right)
(y)\right|  =\mathrm{O}\left(  \varepsilon^{-m\left(  \alpha\right)  }\right)
$ for $\varepsilon\rightarrow0$.

Since $\partial^{\beta}\left(  f\ast\partial^{\alpha}\theta_{\varepsilon
}\right)  =f\ast\partial^{\alpha+\beta}\theta_{\varepsilon}$ the same
arguments applies to derivatives\ and the claim follows.\medskip

Let us now prove that $\iota$ is injective, \textit{id est}
\[
\left(  T\ast\theta_{\varepsilon}\right)  _{\varepsilon}\in\mathcal{N}\left(
\mathrm{C}^{\infty}\left(  \mathbb{R}^{d}\right)  \right)  \Rightarrow
T=0\text{.}%
\]
Indeed, taking $\varphi\in\mathcal{D}\left(  \mathbb{R}^{d}\right)  $ we have
$\left\langle T\ast\theta_{\varepsilon},\varphi\right\rangle \rightarrow
\left\langle T,\varphi\right\rangle $ since $T\ast\theta_{\varepsilon
}\rightarrow T$. But, $T\ast\theta_{\varepsilon}\rightarrow0$ uniformly on
$\operatorname*{supp}\varphi$ since $T\ast\theta_{\varepsilon}\in
\mathcal{N}\left(  \mathrm{C}^{\infty}\left(  \mathbb{R}^{d}\right)  \right)
$.\ Then $\left\langle T\ast\theta_{\varepsilon},\varphi\right\rangle
\rightarrow0$ and $\left\langle T,\varphi\right\rangle =0$.\medskip

We shall prove the last assertion in the case $d=1$, the general case only
differs by more complicate algebraic expressions.

Let $f$ be in $\mathrm{C}^{\infty}\left(  \mathbb{R}\right)  $ and set
$\Delta=\iota_{S}\left(  f\right)  -\sigma\left(  f\right)  .$ One
representative of $\Delta$ is given by
\[
\Delta_{\varepsilon}:\mathbb{R}\rightarrow\mathcal{F}\left(  \mathrm{C}%
^{\infty}\left(  \mathbb{R}\right)  \right)  \;\;\;\;\;y\rightarrow\left(
f\ast\theta_{\varepsilon}\right)  \left(  y\right)  -f(y)=\int f(y-x)\theta
_{\varepsilon}(x)\,\mathrm{d}x-f(y).
\]
Fix $K$ a compact of $\mathbb{R}$. Writing $\int\theta_{\varepsilon}\left(
x\right)  \,\mathrm{d}x=1+\mathcal{N}_{\varepsilon}$ with $\left(
\mathcal{N}_{\varepsilon}\right)  \in\mathcal{N}\left(  \mathbb{R}\right)  $
we get
\[
\Delta_{\varepsilon}(y)=\int\left(  f(y-x)-f(y)\right)  \theta_{\varepsilon
}(x)\,\mathrm{d}x+\mathcal{N}_{\varepsilon}f(y).
\]
The integration is performed on the compact set $\operatorname*{supp}%
\theta_{\varepsilon}\subset\left[  -2/\left|  \ln\varepsilon\right|
,2/\left|  \ln\varepsilon\right|  \right]  $.

Let $k$ be an integer. Taylor's formula gives
\[
f(y-x)-f(y)=\sum_{i=1}^{k}\frac{\left(  -x\right)  ^{i}}{i!}f^{\left(
i\right)  }\left(  y\right)  +\frac{\left(  -x\right)  ^{k}}{k!}\int_{0}%
^{1}f^{\left(  k+1\right)  }\left(  y-ux\right)  \left(  1-u\right)
^{k}\,\mathrm{d}u
\]
and
\begin{multline*}
\Delta_{\varepsilon}(y)=\underset{P_{\varepsilon}\left(  k,y\right)
}{\underbrace{\sum_{i=1}^{k}\frac{\left(  -1\right)  ^{i}}{i!}f^{\left(
i\right)  }\left(  y\right)  \int_{-2/\left|  \ln\varepsilon\right|
}^{2/\left|  \ln\varepsilon\right|  }x^{i}\theta_{\varepsilon}(x)\,\mathrm{d}%
x}}\\
+\underset{R_{\varepsilon}(k,y)}{\underbrace{\int_{-2/\left|  \ln
\varepsilon\right|  }^{2/\left|  \ln\varepsilon\right|  }\frac{\left(
-x\right)  ^{k}}{k!}\int_{0}^{1}f^{\left(  k+1\right)  }\left(  y-ux\right)
\left(  1-u\right)  ^{k}\,\mathrm{d}u\,\theta_{\varepsilon}(x)\,\mathrm{d}x}%
}+\mathcal{N}_{\varepsilon}f(y).
\end{multline*}
According to lemma \ref{SLmnGoodM}, we have $\left(  \int x^{i}\theta
_{\varepsilon}(x)\,\mathrm{d}x\right)  _{\varepsilon}\in\mathcal{N}\left(
\mathbb{R}\right)  $ and consequently
\[
\left(  P_{\varepsilon}\left(  k,y\right)  \right)  _{\varepsilon}%
\in\mathcal{N}\left(  \mathbb{R}\right)  .
\]
Using the definition of $\theta_{\varepsilon}$, we have
\[
R_{\varepsilon}(k,y)=\frac{1}{\varepsilon}\int_{-2/\left|  \ln\varepsilon
\right|  }^{2/\left|  \ln\varepsilon\right|  }\frac{\left(  -x\right)  ^{k}%
}{k!}\int_{0}^{1}f^{\left(  k+1\right)  }\left(  y-ux\right)  \left(
1-u\right)  ^{k}\,\mathrm{d}u\,\rho\left(  \frac{x}{\varepsilon}\right)
\chi\left(  x\left|  \ln\varepsilon\right|  \right)  \,\mathrm{d}x.
\]
Setting $v=x/\varepsilon$ we get
\[
R_{\varepsilon}(k,y)=\varepsilon^{k+1}\int_{-2/\left(  \varepsilon\left|
\ln\varepsilon\right|  \right)  }^{2/\left(  \varepsilon\left|  \ln
\varepsilon\right|  \right)  }\frac{\left(  -v\right)  ^{k}}{k!}\int_{0}%
^{1}f^{\left(  k+1\right)  }\left(  y-\varepsilon uv\right)  \left(
1-u\right)  ^{k}\,\mathrm{d}u\,\rho\left(  v\right)  \chi\left(
\varepsilon\left|  \ln\varepsilon\right|  v\right)  \,\mathrm{d}v.
\]
For $\left(  u,v\right)  \in\left[  0,1\right]  \times\left[  -2/\left(
\varepsilon\left|  \ln\varepsilon\right|  \right)  ,2/\left(  \varepsilon
\left|  \ln\varepsilon\right|  \right)  \right]  $, we have $y-\varepsilon
uv\in\left[  y-1,y+1\right]  $ for $\varepsilon$ small enough. Then, for $y\in
K$, $y-\varepsilon uv$ lies in a compact $K^{\prime}$ for $\left(  u,v\right)
$ in the domain of integration.

It follows
\begin{align*}
\left|  R_{\varepsilon}(k,y)\right|   &  \leq\frac{\varepsilon^{k}}{k!}%
\sup_{\xi\in K^{\prime}}\left|  f^{\left(  k+1\right)  }\left(  \xi\right)
\right|  \int_{-2/\left(  \varepsilon\left|  \ln\varepsilon\right|  \right)
}^{2/\left(  \varepsilon\left|  \ln\varepsilon\right|  \right)  }\left|
v\right|  ^{k+1}\left|  \rho\left(  v\right)  \right|  \mathrm{d}v,\\
&  \leq\frac{\varepsilon^{k}}{k!}\sup_{\xi\in K^{\prime}}\left|  f^{\left(
k+1\right)  }\left(  \xi\right)  \right|  \int_{-\infty}^{+\infty}\left|
v\right|  ^{k+1}\left|  \rho\left(  v\right)  \right|  \mathrm{d}v\leq
C\varepsilon^{k}\;\;(C>0).
\end{align*}
The constant $C$ depends only on the integer $k$, the compacts $K$,
$K^{\prime}$, $\rho$ and $f$.

Finally, for all $k>0$%
\[
\sup_{y\in K}\Delta_{\varepsilon}(y)=\mathrm{o}\left(  \varepsilon^{k}\right)
\;\text{for}\;\varepsilon\rightarrow0.
\]

As $\left(  \Delta_{\varepsilon}\right)  _{\varepsilon}\in\mathcal{F}\left(
\mathrm{C}^{\infty}\left(  \mathbb{R}^{d}\right)  \right)  $ and $\sup_{y\in
K}\Delta_{\varepsilon}(y)=\mathrm{o}\left(  \varepsilon^{k}\right)  $ for all
$k>0$ and $K\Subset\mathbb{R}$, we can conclude directly (without estimating
the derivatives) that $\left(  \Delta_{\varepsilon}\right)  _{\varepsilon}%
\in\mathcal{N}\left(  \mathrm{C}^{\infty}\left(  \mathbb{R}^{d}\right)
\right)  $ by using theorem 1.2.3. of \cite{GKOS}.
\end{proof}

\section{Comparison between $i_{A}$ and $i_{S}$\label{SPComp}}

As quoted in the introduction, the embedding $\iota_{A}:\mathcal{D}^{\prime
}\left(  \mathbb{R}^{d}\right)  \rightarrow\mathcal{G}\left(  \mathbb{R}%
^{d}\right)  $ constructed in \cite{GKOS} depends on the choice of the chosen
net $\rho\in\mathcal{S}\left(  \mathbb{R}^{d}\right)  $. This dependance is a
well known fact for the simplified Colombeau algebra. Of course, $\iota_{S}$
depends also on the choice of $\rho\in\mathcal{S}\left(  \mathbb{R}%
^{d}\right)  $, but not on the choice of $\chi$.\ Moreover:

\begin{proposition}
\label{Sequality}For the same choice of $\rho$, we have: $\iota_{A}=\iota_{S}
$.
\end{proposition}

The proof is carried out in the two following subsections.

\subsection{The case of the embedding of $\mathcal{E}^{\prime}\left(
\mathbb{R}^{d}\right)  $ in $\mathcal{G}\left(  \mathbb{R}^{d}\right)  $}

In \cite{GKOS}, the embedding of $\mathcal{E}^{\prime}\left(  \mathbb{R}%
^{d}\right)  $ in $\mathcal{G}\left(  \mathbb{R}^{d}\right)  $ is realized
with the map
\[
\iota_{0}:\mathcal{E}^{\prime}\left(  \mathbb{R}^{d}\right)  \rightarrow
\mathcal{G}\left(  \mathbb{R}^{d}\right)  \;\;\;\;\;T\mapsto\left(  T\ast
\rho_{\varepsilon}\right)  _{\varepsilon}+\mathcal{N}\left(  \mathrm{C}%
^{\infty}\left(  \mathbb{R}^{d}\right)  \right)  .
\]
We compare here with $\iota_{S\left|  \mathcal{E}^{\prime}\left(
\mathbb{R}^{d}\right)  \right.  }$. Let fix $T\in\mathcal{E}^{\prime}\left(
\mathbb{R}^{d}\right)  $. We have to estimate $\left(  T\ast\rho_{\varepsilon
}\right)  _{\varepsilon}-\left(  T\ast\theta_{\varepsilon}\right)
_{\varepsilon}$.\medskip

Using theorem \ref{SlocStrucEP}, we can write $T=\sum_{finite}\partial
^{\alpha}f_{\alpha}$, each $f_{\alpha}$ having a compact support. Using
linearity, we only need to estimate for one summand, and we shall consider
that $T=\partial^{\alpha}f$. Setting $\Delta_{\varepsilon}=T\ast
\theta_{\varepsilon}-T\ast\rho_{\varepsilon}$ we have
\[
\forall y\in\mathbb{R}^{d},\;\;\;\;\Delta_{\varepsilon}\left(  y\right)  =\int
f(y-x)\left(  \partial^{\alpha}\theta_{\varepsilon}\left(  x\right)
-\partial^{\alpha}\rho_{\varepsilon}\left(  x\right)  \right)  \,\mathrm{d}x.
\]
Then
\begin{align*}
\left|  \Delta_{\varepsilon}\left(  y\right)  \right|   &  \leq C\int\left|
\partial^{\alpha}\theta_{\varepsilon}\left(  x\right)  -\partial^{\alpha}%
\rho_{\varepsilon}\left(  x\right)  \right|  \,\mathrm{d}%
x\;\;\;\;\;\text{\ with }C=\sup_{\xi\in\mathbb{R}^{d}}\left|  f\left(
\xi\right)  \right|  ,\\
&  \leq C\int_{\mathbb{R}^{d}\backslash B(0,1/\left|  \ln\varepsilon\right|
)}\left|  \partial^{\alpha}\theta_{\varepsilon}\left(  x\right)
-\partial^{\alpha}\rho_{\varepsilon}\left(  x\right)  \right|  \,\mathrm{d}x.
\end{align*}
since $\partial^{\alpha}\theta_{\varepsilon}=\partial^{\alpha}\rho
_{\varepsilon}$ on $B(0,1/\left|  \ln\varepsilon\right|  )$.\smallskip

To simplify notations, we suppose $d=1$ and for example $\alpha=1$. We have
\[
\theta_{\varepsilon}^{\prime}\left(  x\right)  -\rho_{\varepsilon}^{\prime
}\left(  x\right)  =\varepsilon^{-1}\left|  \ln\varepsilon\right|
\chi^{\prime}\left(  x\left|  \ln\varepsilon\right|  \right)  \rho\left(
\varepsilon^{-1}x\right)  +\varepsilon^{-2}\rho^{\prime}\left(  \varepsilon
^{-1}x\right)  \left(  \chi\left(  x\left|  \ln\varepsilon\right|  \right)
-1\right)  .
\]
Since $\rho\in\mathcal{S}\left(  \mathbb{R}\right)  $, for all $k\in\mathbb{N}
$, with $k\geq2$, their exists $C(k)\in\mathbb{R}_{+}$ such that
\[
\left|  \rho^{(i)}\left(  x\right)  \right|  \leq\frac{C(k)}{1+\left|
x\right|  ^{k}}\;\;\text{\ (for }i=0\text{ and }i=1\text{).}%
\]
Then, for all $x$ with $\left|  x\right|  \geq1/\left|  \ln\varepsilon\right|
$ we get
\[
\left|  \rho^{(i)}\left(  \varepsilon^{-1}x\right)  \right|  \leq
C(k)\varepsilon^{k}\frac{1}{\varepsilon^{k}+\left|  x\right|  ^{k}}\leq
C(k)\varepsilon^{k}\left|  x\right|  ^{-k}.
\]
Since $\left|  \ln\varepsilon\right|  \leq\varepsilon^{-1}$ for $\varepsilon
\in\left(  0,1\right]  $ and $\left|  \chi\left(  x\left|  \ln\varepsilon
\right|  \right)  -1\right|  \leq1$ for all $x\in\mathbb{R}$, we get
\begin{align*}
\left|  \theta_{\varepsilon}^{\prime}\left(  x\right)  -\rho_{\varepsilon
}^{\prime}\left(  x\right)  \right|   &  \leq\left|  x\right|  ^{-k}\left(
\varepsilon^{k-1}\left|  \ln\varepsilon\right|  \sup\nolimits_{\xi
\in\mathbb{R}}\left|  \chi^{\prime}\left(  \xi\right)  \right|  \left|
x\right|  ^{-k}C(k)\left|  x\right|  ^{-k}+\varepsilon^{k-2}C(k)\left|
x\right|  ^{-k}\right)  ,\\
&  \leq\varepsilon^{k-2}C(k)\left(  \sup\nolimits_{\xi\in\mathbb{R}}\left|
\chi^{\prime}\left(  \xi\right)  \right|  +1\right)  \left|  x\right|  ^{-k}.
\end{align*}
Then, we get a constant $C^{\prime}=C^{\prime}\left(  k,\chi,f\right)  >0$
such that
\[
\left|  \Delta_{\varepsilon}\left(  y\right)  \right|  \leq2\varepsilon
^{k-2}C^{\prime}\int_{1/\left|  \ln\varepsilon\right|  }^{+\infty}\left|
x\right|  ^{-k}\,\mathrm{d}x=\frac{2C^{\prime}}{k-1}\varepsilon^{k-2}\left|
\ln\varepsilon\right|  ^{k-1}.
\]
Finally, we have $\sup_{y\in\mathbb{R}}\left|  \Delta_{\varepsilon}\left(
y\right)  \right|  =\mathrm{o}\left(  \varepsilon^{k}\right)  $ for all
$k\in\mathbb{N}$.

As $\left(  \Delta_{\varepsilon}\right)  _{\varepsilon}\in\mathcal{F}\left(
\mathrm{C}^{\infty}\left(  \mathbb{R}^{d}\right)  \right)  $, we finally
conclude that $\left(  \Delta_{\varepsilon}\right)  _{\varepsilon}%
\in\mathcal{N}\left(  \mathrm{C}^{\infty}\left(  \mathbb{R}^{d}\right)
\right)  $ by using theorem 1.2.3. of \cite{GKOS}. Then:

\begin{lemma}
\label{SLequalityEP}For the same choice of $\rho$, we have: $\iota_{0}%
=\iota_{S\left|  \mathcal{E}^{\prime}\left(  \mathbb{R}^{d}\right)  \right.
}$.
\end{lemma}

\subsection{The case of the embedding of $\mathcal{D}^{\prime}\left(
\mathbb{R}^{d}\right)  $ in $\mathcal{G}\left(  \mathbb{R}^{d}\right)
$\label{SBembRd}}

\begin{notation}
In this subsection we shall note $\mathbb{N}_{m}=\left\{  1,\ldots,m\right\}
$ for all $m\in\mathbb{N}\backslash\left\{  0\right\}  .$
\end{notation}

Let us recall briefly the construction of \cite{GKOS}. Fix some locally finite
open covering $\left(  \Omega_{\lambda}\right)  _{\lambda\in\Lambda}$ with
$\overline{\Omega}_{\lambda}\Subset\mathbb{R}^{d}$ and a family $\left(
\psi_{\lambda}\right)  _{\lambda\in\Lambda}\in\mathcal{D}\left(
\mathbb{R}^{d}\right)  ^{\Lambda}$ with $0\leq\psi_{\lambda}\leq1$ and
$\psi_{\lambda}\equiv1$ on a neighborhood of $\overline{\Omega}_{\lambda} $.
For each $\lambda$ define
\[
\iota_{\lambda}:\mathcal{D}^{\prime}\left(  \mathbb{R}^{d}\right)
\rightarrow\mathcal{G}\left(  \Omega_{\lambda}\right)  \;\;\;\;T\rightarrow
\iota_{\lambda}\left(  T\right)  =\iota_{0}\left(  \psi_{\lambda}T\right)
_{\left|  \Omega_{\lambda}\right.  }=\left(  \left(  \psi_{\lambda}T\ast
\rho_{\varepsilon}\right)  _{\left|  \Omega_{\lambda}\right.  }\right)
_{\varepsilon}+\mathcal{N}\left(  \mathrm{C}^{\infty}\left(  \Omega_{\lambda
}\right)  \right)  .
\]
The family $\left(  \iota_{\lambda}\right)  _{\lambda\in\Lambda}$ is coherent
and by sheaf argument, there exists a unique $\iota_{A}:\mathcal{D}^{\prime
}\left(  \mathbb{R}^{d}\right)  \rightarrow\mathcal{G}\left(  \mathbb{R}%
^{d}\right)  $ such that
\[
\forall\lambda\in\Lambda,\;\;\;\;\iota_{A\left|  \Omega_{\lambda}\right.
}=\iota_{\lambda}.
\]
Moreover, an explicit expression of $\iota_{A}$ can be given: Let $\left(
\chi_{j}\right)  _{j\in\mathbb{N}}$ be a smooth partition of unity subordinate
to $\left(  \Omega_{\lambda}\right)  _{\lambda\in\Lambda}$.\ We have
\[
\forall T\in\mathcal{D}^{\prime}\left(  \mathbb{R}^{d}\right)  ,\;\;\;\;\iota
_{A}\left(  T\right)  =\left(  \sum\nolimits_{j=1}^{+\infty}\chi_{j}\left(
\left(  \psi_{\lambda(j)}T\right)  \ast\rho_{\varepsilon}\right)  \right)
_{\varepsilon}+\mathcal{N}\left(  \mathrm{C}^{\infty}\left(  \mathbb{R}%
^{d}\right)  \right)  .
\]

Let us compare $\iota_{A}$ and $\iota_{S}$.\ Using sheaf properties, we only
need to verify that
\[
\forall\lambda\in\Lambda,\;\;\iota_{s\,\left|  \Omega_{\lambda}\right.
}=\iota_{A\,\left|  \Omega_{\lambda}\right.  }\,\,\,\left(  =\iota_{\lambda
}\right)
\]

For a fixed $\lambda\in\Lambda$ and $T\in\mathcal{D}^{\prime}\left(
\mathbb{R}^{d}\right)  $, we have $\iota_{\lambda}\left(  T\right)  =\iota
_{0}\left(  \psi_{\lambda}T\right)  _{\left|  \Omega_{\lambda}\right.  }$ and
\[
\iota_{A\,\left|  \Omega_{\lambda}\right.  }-\iota_{s\,\left|  \Omega
_{\lambda}\right.  }\left(  T\right)  =\iota_{0}\left(  \psi_{\lambda
}T\right)  -\iota_{s}\left(  \psi_{\lambda}T\right)  +\iota_{s}\left(
\psi_{\lambda}T\right)  -\iota_{s}\left(  T\right)  .
\]
(We omit the restriction symbol in the right hand side).\medskip

As $\psi_{\lambda}T\in\mathcal{E}^{\prime}\left(  \mathbb{R}^{d}\right)  $, we
have $\iota_{0}\left(  \psi_{\lambda}T\right)  =\iota_{s}\left(  \psi
_{\lambda}T\right)  $ according to lemma \ref{SLequalityEP}.\ It remains to
show that $\iota_{s}\left(  \psi_{\lambda}T\right)  =\iota_{s}\left(
T\right)  $, that is to compare $\left(  \left(  \psi_{\lambda}T\right)
\ast\theta_{\varepsilon}\right)  _{\varepsilon}$ and $\left(  T\ast
\theta_{\varepsilon}\right)  _{\varepsilon}$. Let us recall that
\[
\forall y\in\Omega_{\lambda},\;\;\left(  \left(  \psi_{\lambda}T\right)
\ast\theta_{\varepsilon}\right)  _{\varepsilon}(y)-\left(  T\ast
\theta_{\varepsilon}\right)  _{\varepsilon}(y)=\left\langle \psi_{\lambda
}T-T,\left\{  x\mapsto\theta_{\varepsilon}\left(  y-x\right)  \right\}
\right\rangle ,
\]
for $\varepsilon$ small enough.

Let consider $K$ a compact included in $\Omega_{\lambda}$. According to
relation (\ref{SsuppTheta}) we have $\operatorname*{supp}\theta_{\varepsilon
}\left(  y-\cdot\right)  \subset B(y,2/\left|  \ln\varepsilon\right|
)$.\ Using the fact that $\Omega_{\lambda}$ is open, we obtain that
\[
\forall y\in K,\;\;\;\exists\varepsilon_{y}\in\left(  0,1\right]
,\;\;\;\forall\varepsilon\in\left(  0,\varepsilon_{y}\right]
,\;\;\;B(y,2/\left|  \ln\varepsilon\right|  )\subset\Omega_{\lambda}.
\]
The family $\left(  B(y,1/\left|  \ln\varepsilon_{y}\right|  )\right)  _{y\in
K}$ is an open covering of $K$ from which we can extract a finite one,
$\left(  B(y_{l},1/\left|  \ln\varepsilon_{l}\right|  )\right)  _{1\leq l\leq
n} $ (with $\varepsilon_{l}=\varepsilon_{y_{l}}$). Put
\[
\varepsilon_{K}=\min_{1\leq l\leq n}\varepsilon_{l}.
\]
For $y\in K$, there exists $l\in\mathbb{N}_{n}$ such that $y\in B(y_{l}%
,1/\left|  \ln\varepsilon_{l}\right|  )$ Then, for $\varepsilon\leq
\varepsilon_{K}^{2}$, we have$\;$
\[
\operatorname*{supp}\theta_{\varepsilon}\left(  y-\cdot\right)  \subset
B(y,2/\left|  \ln\varepsilon\right|  )\subset B(y,1/\left|  \ln\varepsilon
_{K}\right|  )\subset B(y_{l},2/\left|  \ln\varepsilon_{l}\right|
)\subset\Omega_{\lambda}.
\]
since $d\left(  y,y_{l}\right)  <1/\left|  \ln\varepsilon_{l}\right|  $.

For all $y\in K$, $\theta_{\varepsilon}\left(  y-\cdot\right)  \in
\mathcal{D}\left(  \Omega_{\lambda}\right)  $ for $\varepsilon\in\left(
0,\varepsilon_{K}^{2}\right]  $. Since $T_{\left|  \Omega_{\lambda}\right.
}=\left(  \psi_{\lambda}T\right)  _{\left|  \Omega_{\lambda}\right.  }$ we
finally obtain
\[
\forall y\in K,\;\forall\varepsilon\in\left(  0,\varepsilon_{K}^{2}\right]
,\;\;\left\langle \psi_{\lambda}T-T,\left\{  x\mapsto\theta_{\varepsilon
}\left(  y-x\right)  \right\}  \right\rangle =0,
\]
this showing that $\left(  \left(  \psi_{\lambda}T-T\right)  \ast
\theta_{\varepsilon}\right)  $ lies in $\mathcal{N}\left(  \mathrm{C}^{\infty
}\left(  \mathbb{R}^{d}\right)  \right)  $.

\section{Embedding of $\mathcal{D}^{\prime}\left(  \Omega\right)  $ into
$\mathcal{G}\left(  \Omega\right)  $\label{SPEmbOm}}

All the embeddings of $\mathcal{D}^{\prime}\left(  \Omega\right)  $ into
$\mathcal{G}\left(  \Omega\right)  $ considered in the literature are based on
convolution of distributions by \textrm{C}$^{\infty}$ functions. This product
is possible under additional assumptions in particular about supports. Let
consider both constructions compared in this paper.

For the construction of \cite{GKOS}, the local construction with cutoff
technics applied to the elements of $\mathcal{D}^{\prime}\left(
\Omega\right)  $ is needed to obtain a well defined product of convolution
between elements of $\mathcal{E}^{\prime}\left(  \mathbb{R}^{d}\right)  $ and
$\mathcal{S}\left(  \mathbb{R}^{d}\right)  $. Note that the cutoff is fixed
once for all, and in particular does not depend on $\varepsilon$.

The construction of \cite{NePiSc} allows a ``global'' embedding of
$\mathcal{D}^{\prime}\left(  \mathbb{R}^{d}\right)  $ into $\mathcal{G}\left(
\mathbb{R}^{d}\right)  $ since the convolution of elements of $\mathcal{D}%
^{\prime}\left(  \mathbb{R}^{d}\right)  $ with $\left(  \theta_{\varepsilon
}\right)  _{\varepsilon}\in\left(  \mathcal{D}\left(  \mathbb{R}^{d}\right)
\right)  ^{\left(  0,1\right]  }$ is well defined. But, for the case of an
open subset $\Omega\varsubsetneq\mathbb{R}^{d}$, previous arguments show that
for $y\in\Omega$, the functions $\left\{  x\rightarrow\theta_{\varepsilon
}\left(  y-x\right)  \right\}  $ belongs to $\mathcal{D}\left(  \Omega\right)
$ for $\varepsilon$ smaller than some $\varepsilon_{y}$ depending on $y$. This
does not allow the definition of the net $\left(  T\ast\theta_{\varepsilon
}\right)  _{\varepsilon}$ for $T\in\mathcal{D}^{\prime}\left(  \Omega\right)
$ not compactly supported.\ To overcome this difficulty, a net of cutoffs
$\left(  \kappa_{\varepsilon}\right)  \in\left(  \mathcal{D}\left(
\mathbb{R}^{d}\right)  \right)  ^{\left(  0,1\right]  }$ such that
$\kappa_{\varepsilon}T\rightarrow T$ in $\mathcal{D}^{\prime}\left(
\Omega\right)  $ is considered, giving a well defined convolution of elements
of $\mathcal{E}^{\prime}\left(  \mathbb{R}^{d}\right)  $ with elements of
$\mathcal{D}\left(  \mathbb{R}^{d}\right)  $. We present this construction
below with small changes and an another construction mixing local technics and
compactly supported mollifiers of \cite{NePiSc}.

\subsection{Embedding using cutoff arguments\label{EOmeGlob}}

Let us fix $\Omega\subset\mathbb{R}^{d}$ an open subset and set, for all
$\varepsilon\in\left(  0,1\right]  $,
\[
K_{\varepsilon}=\left\{  x\in\Omega\,\left|  \,d(x,\mathbb{R}^{d}\text{%
$\backslash$%
}\Omega)\geq\varepsilon\text{ and }d(x,0)\leq1/\varepsilon\right.  \right\}
.
\]
Consider $\left(  \kappa_{\varepsilon}\right)  \in\left(  \mathcal{D}\left(
\mathbb{R}^{d}\right)  \right)  ^{\left(  0,1\right]  }$ such that
\[
\forall\varepsilon\in\left(  0,1\right]  ,\;\;\;0\leq\kappa_{\varepsilon}%
\leq1,\;\;\kappa_{\varepsilon}\equiv1\text{ on }K_{\varepsilon}.
\]
(Such a net $\left(  \kappa_{\varepsilon}\right)  _{\varepsilon}$ is obtained,
for example, by convolution of the characteristic function of $K_{\varepsilon
}$ with a net of mollifiers $\left(  \varphi_{\varepsilon}\right)  \in\left(
\mathcal{D}\left(  \mathbb{R}^{d}\right)  \right)  ^{\left(  0,1\right]  }$
with support decreasing to $\left\{  0\right\}  $.)

\begin{proposition}
\label{SPEMBEDOmega}With notations of lemma \ref{SLmnGoodM}, the map
\begin{equation}
\iota_{S}:\mathcal{D}^{\prime}\left(  \Omega\right)  \rightarrow
\mathcal{G}\left(  \Omega\right)  \;\;\;T\mapsto\left(  \left(  \kappa
_{\varepsilon}T\right)  \ast\theta_{\varepsilon}\right)  _{\varepsilon
}+\mathcal{N}\left(  \mathrm{C}^{\infty}\left(  \Omega\right)  \right)
\label{Spomegadef}%
\end{equation}
is an injective homomorphism of vector spaces.\ Moreover $\iota_{S\left|
\mathrm{C}^{\infty}\left(  \mathbb{R}^{d}\right)  \right.  }=\sigma$.
\end{proposition}

We shall not give a complete proof since it is a slight adaptation of the
proof of proposition \ref{SLmnEmbed}.\ We just quote here the main point. As
seen above, many estimates have to be done on compact sets. Let $K$ be a
compact included in $\Omega$ and $\Omega^{\prime}$ an open set such that
$K\subset\Omega^{\prime}\Subset\Omega$. There exists $\varepsilon_{0}%
\in\left(  0,1\right]  $ such that
\[
\forall\varepsilon\in\left(  0,\varepsilon_{0}\right]  ,\;\;\;\Omega^{\prime
}\subset K_{\varepsilon}.
\]
On one hand this implies that we have $\left(  \kappa_{\varepsilon}T\right)
_{\left|  \Omega^{\prime}\right.  }=\left(  \kappa_{\varepsilon_{0}}T\right)
_{\left|  \Omega^{\prime}\right.  }=T_{\left|  \Omega^{\prime}\right.  }$, for
all $T\in\mathcal{D}^{\prime}\left(  \Omega\right)  $ and $\varepsilon
\in\left(  0,\varepsilon_{0}\right]  $.\ On the other hand, we already noticed
that for $y\in K$, the functions $\left\{  x\rightarrow\theta_{\varepsilon
}\left(  y-x\right)  \right\}  $ belongs to $\mathcal{D}^{\prime}\left(
\Omega^{\prime}\right)  $ for all $\varepsilon\in\left(  0,\varepsilon_{K}%
^{2}\right]  $, $\varepsilon_{K}$ only depending on $K$.

Thus a representative of $\iota_{S}\left(  T\right)  $ is given, for all $y\in
K$, by the convolution of an element of $\mathcal{E}^{\prime}\left(
\mathbb{R}^{d}\right)  $ with an element of $\mathcal{D}\left(  \Omega\right)
$ this being valid for $\varepsilon$ smaller than $\min\left(  \varepsilon
_{0},\varepsilon_{K}^{2}\right)  $ only depending on $K$. Proof of
propositions \ref{SLmnEmbed} and \ref{Sequality} can now be adapted using this remark.

\begin{remark}
\label{SCchoice}For the presentation of the construction of \cite{NePiSc} we
chose to consider first the case $\Omega=\mathbb{R}^{d}$. In fact, we can
unify the construction and consider for all $\Omega$ (included in
$\mathbb{R}^{d}$) the embedding defined by (\ref{Spomegadef}). In the case
$\Omega=\mathbb{R}^{d}$, the cutoff functions $\kappa_{\varepsilon}$ are equal
to one on the closed ball $\overline{B\left(  0,1/\varepsilon\right)  }$.
\end{remark}

\subsection{Embedding using local arguments\label{EOmeLoc}}

Let fix $\Omega$ an open subset of $\mathbb{R}^{d}$. Recall that relation
(\ref{SsuppTheta}) implies that
\[
\forall y\in\Omega,\;\;\;\exists\varepsilon_{y}\in\left(  0,1\right]
,\;\;\;\forall\varepsilon\in\left(  0,\varepsilon_{y}\right]
,\;\;\;\operatorname*{supp}\theta_{\varepsilon}\left(  y-\cdot\right)  \subset
B(y,2/\left|  \ln\varepsilon\right|  )\subset\Omega.
\]
and consequently that $\theta_{\varepsilon}\left(  y-\cdot\right)
\in\mathcal{D}\left(  \Omega\right)  $ for $\varepsilon\in\left(
0,\varepsilon_{y}\right]  $. We consider here a local construction to overcome
the fact that $\varepsilon_{y}$ depends on $y$.

Let $\Omega^{\prime}$ be an open relatively compact subset of $\Omega$. As in
subsection \ref{SBembRd}, we find $\varepsilon_{\Omega^{\prime}}$ such that,
for all $\varepsilon\leq\varepsilon_{\Omega^{\prime}}^{2}$ and $y\in
_{\Omega^{\prime}}$, we have $\operatorname*{supp}\theta_{\varepsilon}\left(
y-\cdot\right)  \subset\Omega$ and $\theta_{\varepsilon}\left(  y-\cdot
\right)  \in\mathcal{D}\left(  \Omega\right)  $. For $T\in\mathcal{D}^{\prime
}\left(  \Omega\right)  $, define, for all $y\in\Omega^{\prime}$
\begin{equation}
T_{\varepsilon}\left(  y\right)  =\left\langle T,\theta_{\varepsilon}\left(
y-\cdot\right)  \right\rangle \text{ for }\varepsilon\in\left(  0,\varepsilon
_{\Omega^{\prime}}^{2}\right]  ,\;\;T_{\varepsilon}\left(  y\right)
=T_{\varepsilon_{\Omega^{\prime}}^{2}}\left(  y\right)  \;\text{ for
}\varepsilon\in\left(  \varepsilon_{\Omega^{\prime}}^{2},1\right]  .
\label{SLmnOm0}%
\end{equation}

\begin{lemma}
\label{SLmnOm1}The map
\[
\iota_{\Omega^{\prime}}:\mathcal{D}^{\prime}\left(  \Omega\right)
\rightarrow\mathcal{G}\left(  \Omega^{\prime}\right)  \;\;\;T\mapsto
T_{\varepsilon}\left(  y\right)  +\mathcal{N}\left(  \mathrm{C}^{\infty
}\left(  \Omega^{\prime}\right)  \right)
\]
is an injective homomorphism of vector spaces.
\end{lemma}

The \textbf{proof} is very similar to proposition \ref{SLmnEmbed}'s one.\medskip

Consider now a locally finite open covering of $\left(  \Omega_{\lambda
}\right)  _{\lambda\in\Lambda}$ with $\overline{\Omega}_{\lambda}\Subset
\Omega$ and set $\iota_{\mathbb{\lambda}}=\iota_{\Omega_{\lambda}}$ for
$\lambda\in\Lambda$.

\begin{lemma}
\label{SLmnOm2}The family $\left(  \iota_{\mathbb{\lambda}}\right)
_{\lambda\in\Lambda}$ is coherent.
\end{lemma}

\begin{proof}
Let us take $\left(  \lambda,\mu\right)  \in\Lambda^{2}$ with $\Omega
_{\lambda}\cap\Omega_{\mu}\neq\varnothing$. We have
\[
\iota_{\mathbb{\lambda\,}\left|  \Omega_{\lambda}\cap\Omega_{\mu}\right.
}=\iota_{\mathbb{\mu\,}\left|  \Omega_{\lambda}\cap\Omega_{\mu}\right.  }%
\]
since, for all $T$ in $\mathcal{D}^{\prime}\left(  \Omega\right)  $
representatives\ of $\iota_{\mathbb{\lambda}}$ and $\iota_{\mathbb{\mu}}$,
written in the form (\ref{SLmnOm0}), are equal for $\varepsilon\leq\min\left(
\varepsilon_{\Omega_{\lambda}}^{2},\varepsilon_{\Omega\mu}^{2}\right)  $.
\end{proof}

By sheaf property of $\mathcal{G}\left(  \Omega\right)  $ there exists a
unique $\iota_{S}^{\prime}:\mathcal{D}^{\prime}\left(  \Omega\right)
\rightarrow\mathcal{G}\left(  \Omega\right)  $ such that $\iota_{S\,\left|
\Omega_{\lambda}\right.  }^{\prime}=\iota_{\mathbb{\lambda}}$ for all
$\lambda\in\Lambda$. Moreover, we can give an explicit formula: If $\left(
\Psi_{\lambda}\right)  _{\lambda\in\Lambda}$ is a partition of unity
subordinate to $\left(  \Omega_{\lambda}\right)  _{\lambda\in\Lambda}$, we
have
\[
\forall T\in\mathcal{D}^{\prime}\left(  \Omega\right)  ,\;\;\;\;\iota
_{S}^{\prime}\left(  T\right)  =%
{\textstyle\sum\nolimits_{\lambda\in\Lambda}}
\Psi_{\lambda}\iota_{\mathbb{\lambda}}\left(  T\right)  .
\]
This map $\iota_{S}^{\prime}$ realize an embedding which does not depend on
the particular choice of $\left(  \Omega_{\lambda}\right)  _{\lambda\in
\Lambda}$ (proof left to the reader).

\begin{remark}
\label{SRemOm}One may think that it is regrettable to come back here to local
arguments, whereas they are avoided with cutoff technic.\ This is partially
true but the advantage of compactly supported mollifiers remains: The
convolution with any distribution is possible.\ This renders the local
arguments very simple.
\end{remark}

\subsection{Final remark}

Let $\Omega$ be an open subset of $\mathbb{R}^{d}$.

\begin{proposition}
\label{SMolequal}For the same choice of $\rho$, we have: $\iota_{A}=\iota
_{S}=\iota_{S}^{\prime}$.
\end{proposition}

With notations of previous sections, we only have to prove the equality on
each open set $\Omega_{\lambda}$, where $\left(  \Omega_{\lambda}\right)
_{\lambda\in\Lambda}$ is a covering of $\Omega$ with relatively compact open
sets. As seen before, we shall have $\left(  \kappa_{\varepsilon}T\right)
_{\left|  \Omega_{\lambda}^{\prime}\right.  }=T_{\left|  \Omega_{\lambda
}^{\prime}\right.  }$ and $\theta_{\varepsilon}\left(  y-\cdot\right)
\in\mathcal{D}\left(  \Omega_{\lambda}^{\prime}\right)  $, for all $y\in
\Omega_{\lambda}$ and $\varepsilon$ small enough. ($\Omega_{\lambda}^{\prime}$
is an open subset relatively compact such that $\Omega_{\lambda}\subset
\Omega_{\lambda}^{\prime}\subset\Omega_{\lambda}$.) This remark lead to our
result, since we obtain for $T\in\mathcal{D}\left(  \Omega^{\prime}\right)  $
representatives for $\iota_{S}(T)$ and $\iota_{S}^{\prime}(T)$ equal for
$\varepsilon$ small enough.

\begin{remark}
\label{SembFinalRemark}~\newline $i.$~Let $\mathcal{B}^{\infty}\left(
\mathbb{R}^{d}\right)  $ is the subset of elements of $\mathcal{S}^{\infty
}\left(  \mathbb{R}^{d}\right)  $ satisfying (\ref{SBinfini}). We saw that
there exists fundamentally one class of embeddings $\left(  \iota_{\rho
}\right)  _{\rho\in\mathcal{B}^{\infty}\left(  \mathbb{R}^{d}\right)  }$ of
$\mathcal{D}^{\prime}\left(  \Omega\right)  $ into $\mathcal{G}\left(
\Omega\right)  $ which renders the diagram \ref{Scomdiag} commutative. For a
fixed $\rho\in\mathcal{B}^{\infty}\left(  \mathbb{R}^{d}\right)  $,
$\iota_{\rho}$ can be described globally using technics of \cite{NePiSc} or
locally using either technics of \cite{GKOS} or of subsection \ref{EOmeLoc} of
this paper. This enlarges the possibilities when questions of embeddings arise
in a mathematical problem.\smallskip\newline $ii.$~As mentioned in the
introduction, $\iota_{0}$ can be considered as an embedding of $\mathcal{E}%
^{\prime}\left(  \mathbb{R}^{d}\right)  $ into $\mathcal{G}_{C}\left(
\mathbb{R}^{d}\right)  $. One has the following commutative diagram
\[%
\begin{tabular}
[c]{ccccc}%
$\mathcal{D}$ & $\overset{c}{\longrightarrow}$ & $\mathcal{E}^{\prime}$ &
$\overset{\iota_{0}=\iota_{S}}{\longrightarrow}$ & $\mathcal{G}_{C}$\\
$\downarrow_{c}$ &  & $\downarrow_{c}$ &  & $\downarrow_{i}$\\
$\mathcal{E}=\mathrm{C}^{\infty}\left(  \Omega\right)  $ & $\overset
{c}{\longrightarrow}$ & $\mathcal{D}^{\prime}$ & $\overset{\iota_{A}=\iota
_{S}=\iota_{S}^{\prime}}{\longrightarrow}$ & $\mathcal{G}$%
\end{tabular}
\]
where $\overset{c}{\longrightarrow}$ denote the classical continuous
embedding, and $i$ the canonical embedding of $\mathcal{G}_{C}$ in
$\mathcal{G}$.
\end{remark}

\end{document}